\documentclass[a4paper]{article}
\usepackage{amsfonts}
\usepackage{amssymb}
\usepackage{amsmath}
\usepackage{latexsym}
\usepackage{graphicx}
\usepackage{amsthm}
\usepackage{calrsfs}
\newtheorem{theorem}{Theorem}[section]
\newtheorem{proposition}[theorem]{Proposition}
\newtheorem{corollary}[theorem]{Corollary}
\newtheorem{lemma}[theorem]{Lemma}
\newtheorem{example}[theorem]{Example}
\newtheorem{problem}[theorem]{Problem}
\newtheorem{definition}[theorem]{Definition}
\newtheorem{remark}[theorem]{Remark}

\begin{document}

\parbox{1mm}

\begin{center}
{\bf {\sc \Large An integral with respect to probabilistic-valued
decomposable measures}}
\end{center}

\vskip 12pt

\begin{center}
{\bf Lenka HAL\v CINOV\'A, Ondrej HUTN\'IK}
\end{center}


\hspace{5mm}\parbox[t]{12cm}{\fontsize{9pt}{0.1in}\selectfont\noindent{\bf
Abstract.} Several concepts of approximate reasoning in
uncertainty processing are linked to the processing of
distribution functions. In this paper we make use of probabilistic
framework of approximate reasoning by proposing a Lebesgue-type
approach to integration of non-negative real-valued functions with
respect to probabilistic-valued decomposable (sub)measures. Basic
properties of the corresponding probabilistic integral are
investigated in detail. It is shown that certain properties, among
them linearity and additivity, depend on the properties of the
underlying triangle function providing (sub)additivity condition
of the considered (sub)measure. It is demonstrated that the
introduced integral brings a new tool in approximate reasoning and
uncertainty processing with possible applications in several
areas.} \vskip 24pt

\section{Introduction}

A numerical (i.e., real-valued) measure is a generalization of the
concept of length, area and volume. Natural properties of these
concepts are described by non-negativity and additivity of
considered set functions. The celebrated Lebesgue measure theory
goes even further assuming a countable additivity or
$\sigma$-additivity, which provides a natural background for
probability theory. However, probability theory based on Lebesgue
measure and integral is too restrictive in many diverse fields of
mathematical, economical and engineering sciences, especially in
approximate reasoning, for instance in several kinds of decision
procedures when considering interactions, etc. Thus, there is a
need for non-additive measure and integral theory.

In practice, many physical measurements can be best modelled by
the concept of numerical measure. However, due to the presence of
noise and error in many measurements it is eligible to replace the
(additive) measures with the set functions assigning to each
measurable subset a distribution function, rather than a
non-negative real number. Thus, the probabilistic point of view
comes again to the game: the importance/diameter/measure of a set
might be represented by a distribution function. This resembles
the original idea of Menger replacing a distance function $d:
\Omega\times \Omega \to \mathbb{R}_+$ with a distribution function
$F_{p,q}: \mathbb{R} \to [0, 1]$ wherein for any number $x$ the
value $F_{p,q}(x)$ describes the probability that the distance
between the points $p$ and $q$ in $\Omega$ is less than $x$.
Recently, probabilistic approaches were successfully applied to
various areas in approximate reasoning as well, e.g., to modelling
uncertain logical arguments~\cite{Hunter}, to  approximations of
incomplete data~\cite{GCK}, or to theory of rough sets~\cite{LY}.
Furthermore, a closely related concept can be found in Moore's
interval mathematics~\cite{MB}, where the use of intervals in data
processing is due to measurement inaccuracy and due to rounding.
Observe that intervals can be considered in distribution function
form linked to random variables uniformly distributed over the
relevant intervals. Thus, several concepts of approximate
reasoning in uncertainty processing are linked to the processing
of distribution functions.

The idea behind the above mentioned observations was elaborated
in~\cite{HutMes} leading to the notion of
\textit{$\tau_T$-submeasure}: a set function $\gamma$ defined on a
ring $\Sigma$ of subsets of a non-empty set $\Omega$ taking values
in the set $\Delta^+$ of distribution functions of non-negative
random variables satisfying "initial" condition
$\gamma_\emptyset=\varepsilon_0$ , "antimonotonicity" property
$\gamma_{E}\geq \gamma_F$ whenever $E,F\in\Sigma$ with $E\subseteq
F$, and "subadditivity" property of the form
$$\gamma_{E\cup F}(x+y) \geq T(\gamma_E(x), \gamma_F(y)), \quad
E,F\in\Sigma,\, x,y>0,$$ with $T$ being a left-continuous t-norm.
Here, $\varepsilon_0$ is the distribution function of Dirac random
variable concentrated at point 0. Naturally, the value
$\gamma_E(x)$ describes the probability that the numerical
(sub)measure of a set $E$ is less than $x$. Moreover, a
fuzzy-number-theoretical interpretation is interesting: the value
$\gamma_E$ may be seen as a non-negative LT-fuzzy number, where
$\tau_T(\gamma_E,\gamma_F)$ corresponds to the $T$-sum of fuzzy
numbers $\gamma_E$ and $\gamma_F$ with
\begin{equation}\label{tauT}
\tau_{T}(G,H)(x) := \sup_{u+v=x} T(G(u),H(v)), \quad
G,H\in\Delta^+,
\end{equation} as a special case of triangle functions on $\Delta^+$, i.e., certain semigroups on $\Delta^+$,
see Section~\ref{sectionprem}.

The study of probabilistic-valued set functions continued in
papers~\cite{HalHutMes2} and~\cite{HalHutMes}, where a more
general concept has been proposed. Motivated by the paper of
Shen~\cite{Shen} we have recently introduced in~\cite{HalHutMol} a
general framework of $\tau$-decomposable set functions which
covers all the mentioned approaches and treat them in a unified
way. Indeed, a set function $\gamma: \Sigma\to\Delta^+$ is said to
be a \textit{$\tau$-decomposable measure} if $\gamma_\emptyset =
\varepsilon_0$ and $\gamma_{E\cup F} = \tau(\gamma_E, \gamma_F)$
with $\tau$ being a triangle function on $\Delta^+$. In fact, this
definition resembles the well-known definitions of
$S$-decomposable and $\oplus$-decomposable numerical measures
studied in the classical and pseudo-analysis theory by many
authors. Here, a triangle function $\tau$ is a natural choice for
"aggregation" of $\gamma_E$ and $\gamma_F$ in order to compare
them with $\gamma_{E\cup F}$. For a $\tau$-decomposable measure
$\gamma$ we expect that the value of $\gamma$ at $E\cup F$ is the
same as its value at $F\cup E$ for disjoint sets $E,F\in\Sigma$,
from which follows that
$\tau(\gamma_E,\gamma_F)=\tau(\gamma_F,\gamma_E)$, i.e., $\tau$
has to be commutative. Moreover, from the natural equality
$\gamma_{(E\cup F)\cup G} = \gamma_{E\cup(F\cup G)}$ (as distance
distribution functions) we obtain
$\tau(\tau(\gamma_E,\gamma_F),\gamma_G)=\tau(\gamma_E,\tau(\gamma_F,\gamma_G))$,
i.e., $\tau$ has to be associative. Since
$\gamma_E=\gamma_{E\cup\emptyset}=\tau(\gamma_E,\gamma_\emptyset)=\tau(\gamma_E,\varepsilon_0)$,
then $\varepsilon_0$ has to be a neutral element of $\tau$.
Indeed, $(\Delta^+,\tau)$ has to be a semigroup. Also,
monotonicity of $\tau$ comes into play to provide monotonicity of
the set function $\gamma$. We summarize that the
probabilistic-valued (sub)measures defined and studied
in~\cite{Lpspaces}, \cite{HalHutMes2}, \cite{HalHutMes},
\cite{HutMes} and~\cite{Shen} are only special cases of
$\tau$-decomposable set functions w.r.t. a different choice of
triangle function $\tau$.

A source of motivation for the present paper may be found
in~\cite{Lpspaces} where the authors define a probabilistic-valued
measure, the notion of probabilistic integral of a measurable
function and the corresponding $L^p$-spaces. In doing so their
considerations are related to (in our terminology)
$\tau_M$-decomposable measures defined on a $\sigma$-algebra of
subsets of $\Omega\neq\emptyset$. So, in this paper we follow the
pattern introduced in~\cite{Lpspaces} to define the probabilistic
integral with respect to an arbitrary $\tau$-decomposable
(sub)measure. On one hand, we describe basic properties mentioned
in~\cite{Lpspaces} in detail, on the other hand we investigate
further features of the integral.

The paper is organized as follows. Section~\ref{sectionprem}
contains preliminary notions, such as distance distribution
functions, their topological and algebraic structure, aggregation
functions and triangle functions as well as their relationships
with probabilistic metric spaces. The important facts about
probabilistic-valued set functions with a number of concrete
examples can be found in Section~\ref{secsetfunctions}. In
Section~\ref{secintegral} we define the probabilistic integral
$\int_E f\,\mathrm{d}\gamma$ of a non-negative function
$f:\Omega\to[0,+\infty[$ w.r.t. a $\tau$-decomposable measure
$\gamma:\Sigma\to\Delta^+$ (with $\tau$ being \textit{an arbitrary
distributive triangle function} on $\Delta^+$) on a set
$E\in\Sigma$ using a construction similar to that of the Lebesgue
integral together with a description of its basic properties.
Section~\ref{sec_new} deals with a different characterization of
the integral using sequences of integrals of simple functions
pointwisely converging to the integrand. An interesting feature of
the  integral is its linearity: the integral of (pointwise) sum of
two non-negative measurable functions is a "sum" (in a
probabilistic sense) of their integrals if and only if $\gamma$ is
a $\tau_M$-decomposable measure w.r.t. the triangle function
$\tau_M$ given by~(\ref{tauT}) with $M$ being the minimum t-norm
$M(x,y)=\min\{x,y\}$. This property certifies the importance of
$\tau_M$-decomposable measure and the corresponding integral in
paper~\cite{Lpspaces} when considering probabilistic $L^p$-spaces.
Moreover, we study a $\tau$-decomposable measure induced by the
integral of some fixed function. Furthermore, it is deduced that
the classical integration theory with respect to a positive
measure is naturally included in our general case. The introduced
integral brings a new tool in approximate reasoning and
uncertainty processing with possible applications in several
areas, for example in multicriteria decision support, as it is
demonstrated in Concluding Remarks.

\section{Preliminaries}\label{sectionprem}

In order to make the exposition self-contained, we remind the
reader the basic notions and constructions used in this paper.

\paragraph{Distribution functions and triangle functions} Let $\Delta$ be the family of all distribution functions on the
extended real line $\overline{\mathbb{R}}:=[-\infty,+\infty]$,
i.e., $F: \overline{\mathbb{R}} \to [0,1]$ is non-decreasing, left
continuous on the real line $\mathbb{R}$ with $F(-\infty)=0$ and
$F(+\infty)=1$. The elements of $\Delta$ are partially ordered by
the usual pointwise order $G\leq H$ if and only if $G(x)\leq H(x)$
for all $x\in \overline{\mathbb{R}}$.

A \textit{distance distribution function} is a distribution
function whose support is a subset of
$\overline{\mathbb{R}}_+:=[0,+\infty]$, i.e., a distribution
function $F: \overline{\mathbb{R}}\to [0,1]$ with $F(0) = 0$. The
class of all distance distribution functions will be denoted by
$\Delta^+$.  Distance distribution functions are a proper tool for
measuring distances in probabilistic metric spaces.

A \textit{triangle function} is a function $\tau: \Delta^+ \times
\Delta^+ \to \Delta^+$ which is symmetric, associative,
non-decreasing in each variable and has $\varepsilon_0$ as the
identity, where $\varepsilon_0$ is the distribution function of
Dirac random variable concentrated at point 0. More precisely, for
$a\in[-\infty,+\infty[$ we put
$$\varepsilon_a(x) :=
\begin{cases}
1& \textrm{for}\,\, x>a, \\
0 & \textrm{otherwise}.
\end{cases}$$ The order on $\Delta^+$ induces an order on the set of
triangle functions, i.e., $\tau_1\leq\tau_2$ if and only if
$\tau_1(G,H)(x)\leq\tau_2(G,H)(x)$ for each $x\in\mathbb{R}_+$ and
each $G,H\in\Delta^+$. For more details on triangle functions we
recommend an overview paper~\cite{SamSem}. The most important
triangle functions are those obtained from certain aggregation
functions, especially t-norms.


A \textit{triangular norm}, shortly a t-norm, is a commutative
lattice ordered semi-group on $[0,1]$ with identity 1. The most
important are the minimum t-norm $M(x,y) := \min\{x,y\}$, the
product t-norm $\Pi(x,y) := xy$, the {\L}ukasiewicz t-norm $W(x,y)
:= \max\{x+y-1,0\}$, and the drastic product t-norm
$$D(x,y) :=
\begin{cases}
\min\{x,y\} & \textrm{for}\,\, \max\{x,y\}=1 \\
0 & \textrm{otherwise}.
\end{cases}$$ For more information about t-norms and their properties
we refer to books~\cite{KMP,SS}. Throughout this paper
$\mathcal{T}$ denotes the class of all t-norms.

\paragraph{Binary operations} Let us denote by $\mathcal{L}$ the set of all
binary operations $L$ on $\overline{\mathbb{R}}_+$ such that

\begin{itemize} \item[(i)] $L$ is commutative and associative; \item[(ii)] $L$ is jointly strictly
increasing, i.e., for all $u_1, u_2, v_1,
v_2\in\overline{\mathbb{R}}_+$ with $u_1<u_2$, $v_1<v_2$ holds
$L(u_1,v_1)<L(u_2,v_2)$; \item[(iii)] $L$ is continuous on
$\overline{\mathbb{R}}_+\times\overline{\mathbb{R}}_+$;
\item[(iv)] $L$ has $0$ as its neutral element.
\end{itemize}\noindent Note that $L\in\mathcal{L}$
is a jointly increasing pseudo-addition on
$\overline{\mathbb{R}}_+$ in the sense of~\cite{SM}. The usual
(class of) examples of operations in $\mathcal{L}$ are
\begin{align*}
K_\alpha(x,y) & := (x^\alpha+y^\alpha)^{\frac{1}{\alpha}}, \quad
\alpha>0, \\ K_\infty(x,y) & := \max\{x,y\}.
\end{align*}

\paragraph{Probabilistic metric spaces} Triangular norms and triangle
functions were originally introduced in the context of
probabilistic metric spaces, cf.~\cite{SS}. Recall that a
\textit{probabilistic metric space} (PM-space, for short) is a
non-empty set $\Omega$ together with a family $\mathcal{F}$ of
probability functions $F_{p,q}\in\Delta^+$ satisfying
\begin{align*}
F_{p,q} & = \varepsilon_0\,\,\,\textrm{if and only if}\,\,\,p=q,\\
F_{p,q} & = F_{q,p},
\end{align*} and the "probabilistic analogue" of the triangle
inequality expressed by
\begin{equation}\label{triangleineq}
F_{p,r}\geq \tau(F_{p,q}, F_{q,r})
\end{equation} with $\tau$ being a triangle function on $\Delta^+$, which
holds for all $p,q,r \in \Omega$. The
inequality~(\ref{triangleineq}) depends on a triangle function. In
his original formulation~\cite{Menger} Menger gave as a
generalized triangle inequality the following $$F_{p,r}(x+y)\geq
T(F_{p,q}(x), F_{q,r}(y))\,\,\,\,\textrm{for all}\,\,\,x,y\geq
0,$$ where $T\in\mathcal{T}$ (a left-continuous one). This
corresponds to the triangle function~(\ref{tauT}). Thus, the
triple $(\Omega, \mathcal{F}, \tau_{T})$ is called a
\textit{Menger PM-space} (under $T$). Considering suitable
operations $L$ replacing the standard addition $+$ on
$\overline{\mathbb{R}}_+$ we obtain a larger class of Menger
PM-spaces $(\Omega, \mathcal{F}, \tau_{L,T})$ with triangle
function $\tau_{L,T}$ given by
\begin{equation}\label{tauLT}
\tau_{L,T}(G,H)(x) := \sup_{L(u,v)=x} T(G(u),H(v)), \quad
G,H\in\Delta^+, L\in\mathcal{L}.
\end{equation}

\paragraph{Linear and metric structure on $\Delta^+$} For a distance distribution
function $G$ and a non-negative constant $c\in \mathbb{R}_{+}$ we
define the multiplication of $G$ by a constant $c$ as follows
\begin{equation}\label{operaciestf1}
\left(c\odot G\right)(x):=
\begin{cases}
\varepsilon_0(x), & c=0, \\
G\left(\frac{x}{c}\right), & \textrm{otherwise}.
\end{cases}\end{equation} \noindent Clearly, $c\odot G\in\Delta^+$.
An addition of distance distribution functions may be defined in
the sense of addition via triangular function $\tau$, i.e., we put
\begin{equation}\label{operaciestf2}
(G\oplus_\tau H) (x) := \tau(G, H)(x).
\end{equation} Clearly, $G\oplus_\tau H\in\Delta^+$ for each triangular function $\tau$. We usually
omit a triangle function from the subscript $\oplus_\tau$ and
write just $\oplus$ when no possible confusion may arise. By
associativity of $\tau$ we may introduce the "sum" of $n$
functions $G_1,\dots, G_n \in \Delta^{+}$ as follows
$$\mathop{\bigoplus}_{k=1}^{n} G_k:= \tau\left(G_1,\mathop{\bigoplus}_{k=2}^{n}G_k\right).$$
The following two "laws" will be important for our purposes: for
each $a, b\in\mathbb{R}_{+}$ and each $G\in\Delta^+$ it holds
$$(a\cdot b)\odot G = a\odot(b\odot G) = b\odot(a\odot G).$$ 
A triangle function $\tau$ such that for each $c\in\mathbb{R}_{+}$
and each $G,H\in\Delta^+$ it holds
$$c\odot(G\oplus_\tau H) = (c\odot G)\oplus_\tau(c\odot H)$$ will
be called a \textit{distributive triangle function}. In fact, this
property depends only on $c\in]0,+\infty[$, because for each
triangle function $\tau$ it holds $\tau(\varepsilon_0,
\varepsilon_0) = \varepsilon_0$. The set of all distributive triangle functions on $\Delta^+$ will be denoted by $\mathfrak{D}(\Delta^+)$. 

The set $\Delta^+$ may be endowed with different metrics. We
consider a mapping $d_S: \Delta^+ \times \Delta^+\to [0,1]$ given
by
$$d_S(G,H) = \inf_{h>0} \left\{G(x-h)-h\leq H(x)\leq
G(x+h)+h;\, x\in\left[0,\frac{1}{h}\right)\right\},$$ which is
called the \textit{Sibley metric} (also a modified L\'evy metric)
on~$\Delta^+$, see~\cite{SS}. Immediately, for $G,H\in\Delta^+$
such that $G\leq H$ it holds $d_S(G,\varepsilon_0)\leq
d_S(H,\varepsilon_0)$. A well-known fact is that this metric
metrizes the topology of weak convergence: a sequence
$(G_n)_1^\infty\in\Delta^+$ is said to be \textit{weakly
convergent} to $G\in\Delta^+$, usually written as
$G_n\mathop{\to}\limits^{w\,} G$, if the sequence
$(G_n(x))_1^\infty$ converges to $G(x)$ for every point $x$ of
continuity of $G$. Moreover, $$G_n\mathop{\to}\limits^{w\,}
G\,\,\,\,\, \textrm{if and only if} \,\,\,\,\, d_S(G_n,G)\to 0.$$
Let us recall that the convergence in every point of continuity of
the function $G$ fails to be equivalent to the convergence in any
point of $]0,+\infty[$, see~\cite{SS}.

\begin{remark}\rm\label{reminfimum}
If $\{G_i: i\in I\}$ is a family of functions from $\Delta^+$,
then a pointwise supremum of this family is always a distance
distribution function. On the other hand, the function $G:
\overline{\mathbb{R}}\to[0,1]$ defined as a pointwise infimum
$$G(x)=\inf\{G_i(x): i\in I\}, \, x\in \mathbb{R},$$ is a non-decreasing
function, but it is not necessarily left-continuous on
$\mathbb{R}$. Taking the left-limit
$$\mathfrak{G}(x):=\lim_{x'\nearrow x} G(x) = \sup_{x'<x} G(x), \quad x\in\mathbb{R},$$
the function $\mathfrak{G}$ belongs to $\Delta^+$ and
$\mathfrak{G} = \inf_{i\in I} G_i$ is the infimum of the family
$\{G_i\}$ in the ordered set $(\Delta^+,\leq)$, see~\cite{Cobzas}.
\end{remark}

Finally, a triangle function is \textit{continuous}, if it is
continuous in the metric space $(\Delta^+, d_S)$. Since
$L\in\mathcal{L}$ and $M$ is a continuous t-norm, then the
operation $\oplus_{\tau_{L,M}}$ is continuous, see
e.g.~\cite[Theorem 7.13]{SamSem}.

\section{Probabilistic-valued decomposable set functions w.r.t. a triangle
function: basic facts and examples}\label{secsetfunctions}

Now we introduce the basic notion of probabilistic decomposable
(sub)measure in its general form. For better readability we also
use the following conventions:
\begin{itemize}
\item[(i)] for a probabilistic-valued set function $\gamma: \Sigma
\to \Delta^+$ we write $\gamma_E(x)$ instead of $\gamma(E)(x)$;
\item[(ii)] since $\Delta^+$ is the set of all distribution
functions with support $\overline{\mathbb{R}}_+$, we state the
expression for a mapping $\gamma: \Sigma \to \Delta^+$ just for
positive values of $x$. In case $x\leq 0$ we always suppose
$\gamma_{\cdot}(x)=0$.
\end{itemize}

\begin{definition}\rm\label{defLA-sub}
Let $\tau$ be a triangle function on $\Delta^+$ and $\Sigma$ be a
ring of subsets of $\Omega\neq \emptyset$. A mapping $\gamma:
\Sigma \to \Delta^+$ with $\gamma_\emptyset = \varepsilon_0$ is
said to be a $\tau$-\textit{decomposable submeasure}, if
$\gamma_{E\cup F} \geq \tau(\gamma_E, \gamma_F)$ for each disjoint
sets $E,F \in \Sigma$. If in the preceding inequality equality
holds, then $\gamma$ is said to be a $\tau$-\textit{decomposable
measure} on $\Sigma$.
\end{definition}

Characterization and basic observations about $\tau$-decomposable
(sub)measures are summarized in the following proposition, for the
proofs we refer to~\cite{HalHutMol}.

\begin{proposition}\label{proposition1}
Let $\tau$ be a triangle function on $\Delta^+$ and $\Sigma$ be a
ring of subsets of $\Omega\neq \emptyset$. Then
\begin{itemize}
\item[(i)] $\gamma$ is a $\tau$-decomposable measure if and only
if $\tau(\gamma_{E\cup F}, \gamma_{E\cap F}) = \tau(\gamma_E,
\gamma_F)$ for each $E,F\in\Sigma$; \item[(ii)] each
$\tau$-decomposable measure $\gamma$ is "antimonotone" on
$\Sigma$, i.e., $\gamma_E\geq\gamma_F$ whenever $E,F\in\Sigma$
such that $E\subseteq F$; \item[(iii)] if $\gamma$ is a
$\tau$-decomposable antimonotone submeasure, then the inequality
$\gamma_{E\cup F} \geq \tau(\gamma_E, \gamma_F)$ holds for
arbitrary sets $E,F \in \Sigma$.
\end{itemize}
\end{proposition}

The definition of a $\tau$-decomposable measure may be strengthen
by considering $\sigma$-additive $\tau$-decomposable measures.
Indeed, for a triangle function $\tau$ introduce the notation
$$\mathop{\bigoplus}_{n=1}^{\infty} G_n:=
\lim_{n\to\infty}\mathop{\bigoplus}_{k=1}^{n} G_k$$ with
$G_k\in\Delta^+$ for each $k\in\mathbb{N}$. Then $\gamma$ is said
to be a \textit{probabilistic-valued $\sigma$-additive
$\tau$-decomposable measure} on a $\sigma$-ring $\Sigma$, if
$$\gamma_{\mathop{\cup}\limits_{n=1}^\infty
E_n}=\bigoplus\limits_{n=1}^\infty\gamma_{E_n}$$ whenever
$E_n\in\Sigma$ for each $n\in\mathbb{N}$.

\begin{definition}\rm\label{continuity_below}
We say that a probabilistic-valued set function $\gamma$ is
\textit{continuous from below}, if
$\gamma_{E_n}\mathop{\to}\limits^{w\,}\gamma_E$ whenever
$E_n\nearrow E$, i.e., $E_n\subseteq E_{n+1}$ and
$\mathop{\cup}\limits_{n=1}^\infty E_n=E$.
\end{definition}

The role of continuity in approximate reasoning is emphasized e.g.
in paper~\cite{Jenei}. In the whole paper we will consider the
limit operation taken in the weak topology on $\Delta^+$, i.e.,
the limit $\lim\limits_{n\to\infty} G_n = G$ is always understood
in the metric space $(\Delta^+, d_S)$. Then the following result
holds.

\begin{lemma}\label{lemma_continuity}
Let $\tau$ be a triangle function on $\Delta^+$ and $\gamma$ be a
$\sigma$-additive $\tau$-decomposable measure on a $\sigma$-ring
$\Sigma$. If $E_n\nearrow E$, then $\gamma$ is continuous from
below.
\end{lemma}

\proof For each $n\in\mathbb{N}$ put $F_n=E_n\setminus E_{n-1}$,
where $E_0:=\emptyset$. Clearly,
$\mathop{\cup}\limits_{n=1}^\infty
F_n=\mathop{\cup}\limits_{n=1}^\infty E_n$ and
$\mathop{\cup}\limits_{k=1}^n F_k=E_n$. So,
$\gamma_{E_n}=\gamma_{\mathop{\cup}\limits_{k=1}^n
F_k}=\bigoplus\limits_{k=1}^n\gamma_{F_k}$ and it follows
$$\gamma_{E}=\gamma_{\mathop{\cup}\limits_{n=1}^\infty
E_n}=\gamma_{\mathop{\cup}\limits_{n=1}^\infty
F_n}=\bigoplus\limits_{n=1}^\infty\gamma_{F_n}=\lim\limits_{n\to\infty}\bigoplus\limits_{k=1}^n\gamma_{F_k}=\lim\limits_{n\to\infty}\gamma_{E_n},$$
which completes the proof. \qed

Let us present some examples of probabilistic-valued
(sub)measures. The first one claims that each numerical measure
may be regarded as a probabilistic-valued decomposable measure. It
is, in fact, highly expected, because similar situation appears in
the case of metric spaces and PM-spaces: each metric space
$(\Omega,d)$ is a Menger PM-space $(\Omega,\mathcal{F},\tau_T)$
with $F_{p,q}=\varepsilon_{d(p,q)}$ and a left-continuous t-norm
$T$, cf.~\cite[Remark 9.25]{KMP}. So, probabilistic-valued
measures offer a wider framework than that of the classical
measures and are general enough to cover even wider statistical
situations. The importance of such measures can also be traced in
other areas, e.g., stochastic differential equations.

\begin{example}\rm\label{example MS}
Let $m$ be a finite numerical measure on $\Sigma$. If $m$ is
$L$-decomposable, i.e., $m(E\cup F)=L(m(E),m(F))$ with
$L\in\mathcal{L}$, then for any $\Phi\in\Delta^+$ the set function
$\gamma^\Phi:\Sigma\to\Delta^+$ defined by
$\gamma^\Phi_E:=m(E)\odot \Phi$ with $E\in\Sigma$ is a
$\tau_{L,M}$-decomposable measure on $\Sigma$, where $\tau_{L,M}$
is given by~(\ref{tauLT}) with the minimum t-norm $M$. Indeed, the
additivity property $\gamma^\Phi_{E\cup
F}=\tau_{L,M}\left(\gamma^\Phi_{E},\gamma^\Phi_{F}\right)$ follows
from~\cite[Section 7.7]{SS}: the equality
\begin{equation}\label{operaciestf3}
\tau\left(c_1\odot G, c_2\odot G\right)= L(c_1,c_2)\odot G
\end{equation} holds for each $c_{1}, c_{2}\in\mathbb{R}_{+}$ and each $G\in
\Delta^+$ if and only if $\tau=\tau_{L,M}$. Moreover,
see~\cite{AS}, $$\tau\left(c_1\odot G, c_2\odot
G\right)\geq(\leq)\, L(c_1,c_2)\odot
G\,\,\,\,\Leftrightarrow\,\,\,\,\tau\geq(\leq)\,\tau_{L,M}.$$
These properties will be crucial when investigating linearity of
the introduced probabilistic integral in
Section~\ref{secintegral}.
\end{example}

\begin{figure}
\begin{center}
\includegraphics{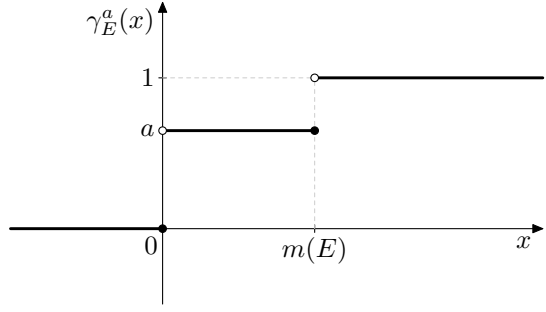}
\caption{Graph of a probabilistic measure $\gamma^a_{E}$ from
Example~\ref{examplegamma}} \label{obrPM}
\end{center}
\end{figure}

\begin{example}\rm\label{examplegamma}
A special case of Example~\ref{example MS} will be useful for us
later to demonstrate certain properties of the introduced integral
in Section~\ref{secintegral}. Indeed, for $a\in[0,1]$ let
$\lambda^a\in \Delta^+$ be defined as follows
\begin{eqnarray}
\lambda^a(x)=\, a \chi_{]0, 1]}(x) + \chi_{]1, +\infty]}(x)\, =\,
\begin{cases}
0, & x\leq 0, \\
a, & x \in ]0, 1] \\ 1, & x> 1,
\end{cases}
\end{eqnarray}
where $\chi_S$ is the characteristic function of a set $S$.
Clearly, $\lambda^0=\varepsilon_1$ and $\lambda^1=\varepsilon_0$.
Moreover, the set $\Lambda=\{\lambda^a; a\in[0,1]\}$ is linearly
ordered, i.e., $\lambda^a\leq\lambda^b$ whenever $a\leq b$. For a
finite numerical (additive) measure $m$ on a ring $\Sigma$ put
$\gamma^a_E:=m(E)\odot\lambda^a$ for $E\in\Sigma$, i.e., we have
\begin{eqnarray}\label{univmiera}
\gamma^a_{E}(x)=\, a \chi_{]0, m(E)]}(x) + \chi_{]m(E),
+\infty]}(x)\, =\,
\begin{cases}
0, & x\leq 0, \\
a, & x \in ]0, m(E)] \\ 1, & x> m(E).
\end{cases}
\end{eqnarray}
Graph of $\gamma^a_{E}$ for some $a\in[0,1]$ and certain value of
measure $m(E)$ is illustrated in Fig.~\ref{obrPM}. Then $\gamma^a$
is a $\tau_{M}$-decomposable measure. For $A \subseteq [0, 1]$ we
denote by $\Gamma_A$ the following set of all such probabilistic
$\tau_M$-decomposable measures.
\end{example}

\begin{example}\rm\label{exShen}
\textit{Shen's $\top$-probabilistic decomposable measures},
cf.~\cite{Shen}: this class of measures
$\mathfrak{M}:\Sigma\to\Delta^+$ of the form $\mathfrak{M}_{E\cup
F}(t)=\top(\mathfrak{M}_{E}(t),\mathfrak{M}_{F}(t))$ for disjoint
$E,F\in\Sigma$ corresponds to the class of $\tau$-decomposable
measures w.r.t. the pointwisely defined function $\Pi_\top:
\Delta^{+}\times \Delta^{+}\rightarrow \Delta^{+}$ of the form
\begin{equation}\label{tautop}
\Pi_\top(G,H)(t)=\top(G(t), H(t)), \quad
G,H\in\Delta^+,\end{equation} with $\top$ being a left-continuous
t-norm. Left-continuity of $\top$ is a necessary and sufficient
condition for $\Pi_\top$ being a triangle function,
cf.~\cite[Theorem 5.2]{SamSem}.
\end{example}

\begin{example}\rm
Let $\mu: \Sigma\to[0,+\infty)$ be a set function with
$\mu(\emptyset)=0$ and for $E\in\Sigma$ put
$\gamma_E=\varepsilon_{\mu(E)}$. Then
\begin{itemize}
\item[(i)] if $\mu$ is additive, then $\gamma$ is a
$\tau_*$-decomposable measure, where
$$\tau_*(G,H)(x):= (G*H)(x)=
\begin{cases}
0, & x=0, \\
\int_{0}^x G(x-t)\,\textrm{d}H(t), & x\in]0,+\infty[, \\
1, & x=+\infty,
\end{cases}$$ is the convolution of $G,H\in\Delta^+$ and the integral is
meant in the sense of Lebesgue-Stieltjes; moreover, $\gamma$ is a
$\tau_D$-decomposable measure (i.e., for non-continuous triangle
function $\tau_D$ with $D$ being the drastic product), and it is a
$\tau_T$-decomposable measure for an arbitrary continuous
$T\in\mathcal{T}$ (i.e., for continuous triangle function
$\tau_T$); \item[(iii)] if $\mu$ is $L$-decomposable with
$L\in\mathcal{L}$, and $T\in\mathcal{T}$ is continuous, then
$\gamma$ is a $\tau_{L,T}$-decomposable measure; \item[(iv)] if
$\mu$ is $K_\infty$-decomposable, then $\gamma$ is
$\Pi_\top$-decomposable measure for each $\top\in\mathcal{T}$,
where $\Pi_\top$ is given by~(\ref{tautop}), i.e., $\gamma$ is a
Shen's $\top$-probabilistic decomposable measure (in the
terminology and notation of~\cite{Shen}).
\end{itemize}
\end{example}


Further examples of probabilistic (sub)measures w.r.t. different
triangle functions may be found in~\cite{HalHutMes2} as well as
in~\cite{HalHutMol}.

\section{Probabilistic integral: a general construction}\label{secintegral}

In what follows let $\Omega^{\overline{\mathbb{R}}_{+}}$ be a set
of all functions $f: \Omega\to\overline{\mathbb{R}}_+$. Motivating
by the concept of "probabilistic integral" introduced
in~\cite{Lpspaces} we aim to investigate here in detail a little
bit more general definition of integral
$$\int_{E} f \,\mathrm{d}\gamma,$$ where

\begin{itemize}
\item $E\in \Sigma$ with $\Sigma$ being a ring of subsets of a
non-empty set $\Omega$;

\item $f\in \Omega^{\overline{\mathbb{R}}_{+}}$ is a measurable
function (in the sense described below);

\item $\gamma: \Sigma \to \Delta^{+}$ is a $\tau$-decomposable
measure with respect to a distributive triangle function $\tau$.
\end{itemize}

In what follows (for the purpose of this section) we always
suppose that \textsf{$\gamma$ is a probabilistic-valued
$\tau$-decomposable measure w.r.t. to an arbitrary distributive
triangle function $\tau$} and all functions are supposed to be
measurable. The introduction of the $\gamma$-integral follows the
pattern known from the classical Lebesgue integral: first it is
defined for measurable step functions and then for non-negative
measurable functions using the fact that each such functions can
be approximated uniformly by a monotone sequence of step
functions.

\begin{definition}\rm
A function $f\in \Omega^{\overline{\mathbb{R}}_{+}}$ is called
\textit{simple}, if it attains a finite number of values $x_1,
x_2,\dots, x_n \in\mathbb{R}_{+}$ on pairwise disjoint sets $E_1,
E_2,\dots, E_n\in\Sigma$. For such a function we write
$f=\sum\limits_{i=1}^{n}x_{i}\chi_{E_{i}}$. A set of all simple
functions from $\Omega^{\overline{\mathbb{R}}_{+}}$ will be
denoted by $\mathcal{S}$.
\end{definition}

Clearly, if $f$ and $g$ are simple functions and
$c\in\mathbb{R}_+$, then $f+g$ and $cf$ are simple functions as
well. The simplest example of a simple function is the
characteristic function of a measurable set $F\in\Sigma$, i.e.,
$$\chi_F(x):=\begin{cases}
1, & x\in F, \\
0, & x\notin F.\end{cases}$$ Its $\gamma$-integral w.r.t. a
$\tau$-decomposable measure $\gamma$ on $\Sigma$ with $\tau$ being
a distributive triangle function is then defined by the equality
$$\int_{E} \chi_F\,\mathrm{d}\gamma := \gamma_{E\cap F}.$$ For $f\in\mathcal{S}$
we put
$$\int_E f\,\mathrm{d}\gamma = \int_E \left(\sum\limits_{i=1}^{n}x_{i}\chi_{E_{i}}\right)\,\mathrm{d}\gamma :=
\mathop{\bigoplus}\limits_{i=1}^{n}x_i\odot\gamma_{E\cap E_{i}}.$$
Clearly, the $\gamma$-integral of a simple function is a distance
distribution function. The correctness of the introduced
definition of integral will be examined in the following lemma
which enlightens why distributivity of a triangle function is
crucial in the construction of the integral.

\begin{lemma}
Let $\tau\in\mathfrak{D}(\Delta^+)$ and $\gamma$ be a
$\tau$-decomposable measure on $\Sigma$. If $f\in\mathcal{S}$ is
of the form $f=\sum\limits_{i=1}^{n}\alpha_{i}\chi_{E_{i}}$, where
$\alpha_{i}\in\mathbb{R}_{+}$ and $E_i\in\Sigma$, $i=1,\dots,n$,
are pairwise disjoint sets, and $g\in\mathcal{S}$ is of the form
$g=\sum\limits_{j=1}^{m}\beta_{j}\chi_{F_{j}}$ where
$\beta_{j}\in\mathbb{R}_{+}$ and $F_j\in\Sigma$, $j=1,\dots,m$,
are pairwise disjoint sets, then
$$\int_E f\,\mathrm{d}\gamma=\int_E g\,\mathrm{d}\gamma$$ provided
that $f=g$.
\end{lemma}

\proof We may consider $\alpha_i\neq 0$ and $\beta_j\neq 0$ for
each $i=1,2,\dots,n$ and $j=1,2,\dots,m$, because multiplication
of a distance distribution function by zero produces the neutral
element $\varepsilon_0$, see~(\ref{operaciestf1}). Since $E_i$,
resp. $F_j$ are disjoint, then $E_i=
\mathop{\cup}\limits_{j=1}^m(E_i\cap F_j)$ for each $i=1,\dots,n$
as well as $F_j=\mathop{\cup}\limits_{i=1}^n(F_j\cap E_i)$ for
each $j=1,\dots,m$. Also, if $E_i\cap F_j\neq\emptyset$, then
$\alpha_i=\beta_j$. Hence, $\alpha_i\odot\gamma_{E_i\cap
F_j}=\beta_j\odot\gamma_{E_i\cap F_j}$ (for $E_i\cap
F_j=\emptyset$ the both sides are equal to $\varepsilon_0$).
Finally, from $\tau$-decomposability of $\gamma$ we get
\begin{align*}\int_E f\,\mathrm{d}\gamma & = \bigoplus_{i=1}^n
\alpha_i\odot\gamma_{E\cap E_i}= \bigoplus_{i=1}^n
\alpha_i\odot\gamma_{E\cap\left(\mathop{\cup}\limits_{j=1}^m(E_i\cap
F_j)\right)}= \bigoplus_{i=1}^n \alpha_i\odot\bigoplus_{j=1}^m
\gamma_{E\cap(E_i\cap F_j)}\\ &=\bigoplus_{i=1}^n
\bigoplus_{j=1}^m \alpha_i\odot\gamma_{E\cap(E_i\cap
F_j)}=\bigoplus_{i=1}^n \bigoplus_{j=1}^m
\beta_j\odot\gamma_{E\cap(E_i\cap F_j)}= \bigoplus_{j=1}^m
\bigoplus_{i=1}^n \beta_j\odot\gamma_{E\cap(E_i\cap
F_j)}\\
&=\bigoplus_{j=1}^m \beta_j \odot \bigoplus_{i=1}^n
\gamma_{E\cap(E_i\cap F_j)}=\bigoplus_{j=1}^m
\beta_j\odot\gamma_{E\cap\left(\mathop{\cup}\limits_{i=1}^n(E_i\cap
F_j)\right)}=\bigoplus_{j=1}^m \beta_j\odot\gamma_{E\cap F_j}
\\ & =\int_E g\,\mathrm{d}\gamma\end{align*} which completes the
proof. \qed \vskip 5 pt

Because of this, the integral of a simple function w.r.t. a
probabilistic-valued $\tau$-decomposable measure $\gamma$ with
$\tau$ being a distributive triangle function is well-defined.

\begin{example}\rm\label{priklad1}
Let $\gamma^a$ be given by~(\ref{univmiera}) for some $a\in[0,1]$.
If $f\in\mathcal{S}$, then the $\gamma^a$-integral of $f$ on
$E\in\Sigma$ is a distribution function of the form
$$\int_E f\,\mathrm{d}\gamma^{a}=\int_E \left(\sum\limits_{i=1}\limits^{n}x_{i}\chi_{E_{i}}\right)\,\mathrm{d}\gamma^{a}=
r\odot\lambda^a,$$ where $r=\sum\limits_{i=1}\limits^{n}
x_{i}\mu(E\cap E_{i})$ with $x_i\in \mathbb{R}_{+}$ and $E_i$,
$i=1,2,\dots, n$, being pairwise disjoint sets from $\Sigma$. The
resulting form of integral follows from addition
$\oplus_{\tau_{M}}$ of distribution functions in Menger PM-space
(under the minimum t-norm $M$), see the
equality~(\ref{operaciestf3}). Note that $\tau_M$ is a
distributive triangle function. Coefficient $r$ may be expressed
in the form $r=\int_E f\,\mathrm{d}\mu$ of the classical Lebesgue
integral. 

More generally, if we consider the operation $\oplus_{\tau}$ with
$\tau=\tau_{L,M}$ given by~(\ref{tauLT}) for $L\in\mathcal{L}$,
the $\gamma^a$-integral of $f$ on $E\in\Sigma$ is the function
$r\odot\lambda^a$ with
$$r= \mathop{\mathcal{L}}_{i=1}^n r_i := L(r_1,
L(r_2,\dots,L(r_{n-1},r_{n})\dots)),\quad r_i=x_{i}\,\mu(E\cap
E_{i}).$$
\end{example}


Let $f$ be a simple non-negative function on $\Omega$. Since
$\int_E 0\,\mathrm{d}\gamma = \varepsilon_0$ and
$H\leq\varepsilon_0$ for each $H\in\Delta^+$, then we immediately
get $\int_E f\,\mathrm{d}\gamma \leq \varepsilon_0$. So, an
interesting feature of the introduced integral of simple functions
is the following "order reversing" property.

\begin{theorem}\label{antimonotonicity}
Let $\tau\in\mathfrak{D}(\Delta^+)$, $\gamma$ be a
$\tau$-decomposable measure on $\Sigma$ and $f, g\in\mathcal{S}$.
If $f\leq g$ then $\int_E f\,\mathrm{d}\gamma\geq \int_E
g\,\mathrm{d}\gamma$ on $E\in\Sigma$.
\end{theorem}

\proof  Let $E\in\Sigma$ and $f,g\in\mathcal{S}$ such that $f\leq
g$. Then there exist pairwise disjoint sets $E_1,\dots,
E_n\in\Sigma$ such that both functions may be expressed in the
form
$$f=\sum\limits_{i=1}^{n}\alpha_{i}\chi_{E_{i}},
\quad g=\sum\limits_{i=1}^{n}\beta_{i}\chi_{E_{i}}$$ with
$\alpha_{i}, \beta_{i}\in\mathbb{R}_+$ and
$0\leq\alpha_i\leq\beta_i<+\infty$. From monotonicity of distance
distribution functions it follows $\alpha_i\odot F\geq
\beta_i\odot F$ for each $F\in\Delta^+$. 
Hence, $$\int_E f\,\mathrm{d}\gamma=\bigoplus_{i=1}^n
\alpha_i\odot\gamma_{E\cap E_i}\geq \bigoplus_{i=1}^n
\beta_i\odot\gamma_{E\cap E_i}=\int_E g\,\mathrm{d}\gamma,$$
completing the proof. \qed \vskip 5 pt

The "order reversing" property of the introduced integral of a
simple function motivates the following definition of the integral
of a non-negative measurable function. Measurability will be
understood in the following sense: a function $f\in
\Omega^{\overline{\mathbb{R}}_{+}}$ is measurable w.r.t. $\Sigma$,
if there exists a sequence $(f_n)_1^\infty\in \mathcal{S}$ such
that $\lim\limits_{n\to\infty} f_n(x)=f(x)$ for each $x\in\Omega$.
For a measurable function $f\in
\Omega^{\overline{\mathbb{R}}_{+}}$ we denote by
$\mathcal{S}_{f,E}$
the set of all simple functions $\mathfrak{f}$ with $\mathfrak{f}(t)\leq f(t)$ for each $t\in E$. 


\begin{definition}\label{integrability0}\rm
Let $\tau\in\mathfrak{D}(\Delta^+)$ and $\gamma$ be a
$\tau$-decomposable measure on $\Sigma$. We say that a function
$f\in \Omega^{\overline{\mathbb{R}}_{+}}$ is
\textit{$\gamma$-integrable} on a set $E\in\Sigma$, if there
exists $H\in\Delta^+$ such that $\int_E
\mathfrak{f}\,\mathrm{d}\gamma\geq H$ for each $\mathfrak{f}\in
\mathcal{S}_{f,E}$. In this case we put
$$\int_E f\,\mathrm{d}\gamma :=
\inf\left\{\int_E \mathfrak{f}\,\mathrm{d}\gamma;\,
\mathfrak{f}\in \mathcal{S}_{f,E}\right\}.$$
\end{definition}

\begin{remark}\rm
We emphasize the fact that the values of a simple function are to
be finite non-negative real numbers, however a $\gamma$-integrable
function is extended non-negative real-valued. Indeed, if
$f(x)=+\infty$, then by Halmos~\cite[\S 20, Theorem B]{Halmos}
$f_n(x)=n$ for every $n$. Thus, putting
$$\varepsilon_\infty(x) :=
\begin{cases}
1& \textrm{for}\,\, x=+\infty, \\
0 & \textrm{otherwise},
\end{cases}$$ and extending the multiplication $\odot$ by
$+\infty$ as follows
$$
\left(c\odot G\right)(x):= \varepsilon_\infty(x), \quad
c=+\infty,$$ the integral of such a function $f$ is obviously
equal to $\varepsilon_\infty$. This enables to integrate extended
non-negative real-valued functions.
\end{remark}

Since all properties of simple (integrable) functions and their
integrals remain valid also for general integrable functions and
their integrals, in theorems we shall use the following
formulation: ...(simple) $\gamma$-integrable functions... On one
hand, by this formulation we want to emphasize that the theorem
has an importance by itself for simple functions, while on the
other hand, we want to emphasize that the theorem is valid for
general integrable functions. We prove them first for simple
functions and at the same time we point out their proofs for
general integrable functions.

\begin{theorem}\label{nvlastnosti}
Let $\tau\in\mathfrak{D}(\Delta^+)$ and $\gamma$ be a
$\tau$-decomposable measure on $\Sigma$. Let $f, g\in
\Omega^{\overline{\mathbb{R}}_{+}}$ be (simple)
$\gamma$-integrable functions on the corresponding sets and put
$\mathcal{I}: f\mapsto\int_{\cdot}{f \,\mathrm{d}\gamma}$. Then
$\gamma$\,-integral is
\begin{itemize}
\item[(i)] an antimonotone operator, i.e., $\mathcal{I}(f)\geq
\mathcal{I}(g)$ whenever $f\leq g$; \item[(ii)] positively
homogeneous, i.e., for each $c\in \mathbb{R}_{+}$ it holds
$\mathcal{I}(c\cdot f) = c \odot \mathcal{I}(f)$.
\end{itemize}
\end{theorem}

\proof (i) For simple functions $f\leq g$ the inequality
$\mathcal{I}(f)\geq \mathcal{I}(g)$ follows from
Theorem~\ref{antimonotonicity}. Let $f,g\in
\Omega^{\overline{\mathbb{R}}_{+}}$ be such that $f\leq g$. Then
$$\left\{\int_E \mathfrak{f}\,\mathrm{d}\gamma;\,\, \mathfrak{f}\in \mathcal{S}_{f,E}\right\}\subseteq \left\{\int_E
\mathfrak{g}\,\mathrm{d}\gamma;\,\, \mathfrak{g}\in
\mathcal{S}_{g,E}\right\},$$ from it follows
$$\inf\left\{\int_E \mathfrak{f}\,\mathrm{d}\gamma;\,\, \mathfrak{f}\in \mathcal{S}_{f,E}\right\}\geq
\inf\left\{\int_E \mathfrak{g}\,\mathrm{d}\gamma;\,\,
\mathfrak{g}\in \mathcal{S}_{g,E}\right\}.$$ This implies the
required inequality.

(ii) Let $c\in \mathbb{R}_{+}$. For a simple function $f$ of the
form $f=\sum\limits_{i=1}^n x_i\chi_{E_i}$ the function $c\cdot f$
is also simple of the form $c\cdot f=\sum\limits_{i=1}^n c
x_i\chi_{E_i}$. Then
$$\int_E c\cdot f\,\mathrm{d}\gamma = \bigoplus_{i=1}^n cx_i\odot\gamma_{E\cap E_i} =
\bigoplus_{i=1}^n c\odot(x_i\odot\gamma_{E\cap E_i}) =
c\odot\bigoplus_{i=1}^n x_i\odot\gamma_{E\cap E_i} = c\odot\int_E
f\,\mathrm{d}\gamma,$$ where the third equality is due to
distributivity of $\tau$. Then for a general $\gamma$-integrable
function $f$ we get
$$\int_E c\cdot f\,\mathrm{d}\gamma = \inf\left\{\int_E
c\cdot \mathfrak{f}\,\mathrm{d}\gamma;\,\, \mathfrak{f}\in
\mathcal{S}_{f,E}\, \right\}= c\odot \inf\left\{\int_E
\mathfrak{f}\,\mathrm{d}\gamma;\,\, \mathfrak{f}\in
\mathcal{S}_{f,E}\right\} = c \odot\int_E f\,\mathrm{d}\gamma,$$
which completes the proof. \qed \vskip 5pt

Now we examine the $\gamma$-integral with respect to the structure
of measures as distance distribution functions: homogeneity and
monotonicity.

\begin{theorem}
Let $\tau\in\mathfrak{D}(\Delta^+)$ and $\gamma^i$ be
$\tau$-decomposable measures on $\Sigma$ with $i\in\{1,2\}$. If
$f\in \Omega^{\overline{\mathbb{R}}_{+}}$ is a (simple)
$\gamma^i$-integrable function on $E\in\Sigma$ for each
$i\in\{1,2\}$, then
\begin{itemize}
\item[(i)] for each $c\in\mathbb{R}_+$ the function $f$ is
$c\odot\gamma^1$-integrable on $E\in\Sigma$, and it holds
$$\int_E f\,\mathrm{d}(c\odot\gamma^1) = c\odot\int_E
f\,\mathrm{d}\gamma^1;$$ \item[(ii)] if
$\gamma^1_E\leq\gamma^2_E$, then
$$\int_E f\,\mathrm{d}\gamma^1 \leq \int_E
f\,\mathrm{d}\gamma^2.$$
\end{itemize}
\end{theorem}

\proof First observe that if $\gamma^i$ are $\tau$-decomposable
measures on $\Sigma$ for $i\in\{1,2\}$, then by~\cite[Theorem
4.2]{HalHutMol} the set function $\gamma:=c\odot\gamma^1$,
$c\in\mathbb{R}_+$ is a $\tau$-decomposable measure on $\Sigma$.

(i) Since for a simple function $f=\sum\limits_{i=1}^n
x_i\chi_{E_i}$ we get
$$\int_E f\,\mathrm{d}\gamma = \bigoplus_{i=1}^n
x_i\odot\left(c\odot\gamma^1_{E\cap E_i}\right) =
\bigoplus_{i=1}^n cx_i\odot\gamma^1_{E\cap E_i} = \int_E c\cdot
f\,\mathrm{d}\gamma^1 = c\odot\int_E f\,\mathrm{d}\gamma^1,$$
where the last equality follows from
Theorem~\ref{nvlastnosti}(ii), for a general non-negative function
$f$ the result holds by Theorem~\ref{nvlastnosti}(ii) as well.

(ii) From monotonicity of operation $\oplus_\tau$ for a simple
function $f=\sum\limits_{i=1}^n x_i\chi_{E_i}$ it holds
\begin{align*}
\int_E f\,\mathrm{d}\gamma^1 &= \bigoplus_{i=1}^n
x_i\odot\gamma^1_{E\cap E_i} \leq \bigoplus_{i=1}^n
x_i\odot\gamma^2_{E\cap E_i}= \int_E
f\,\mathrm{d}\gamma^2.\end{align*} If $f$ is a non-negative
integrable function, then
$$\int_E f\,\mathrm{d}\gamma^1 = \inf\left\{\int_E
\mathfrak{f}\,\mathrm{d}\gamma^1;\,\,
\mathfrak{f}\in\mathcal{S}_{f,E}\right\} \leq \inf\left\{\int_E
\mathfrak{f}\,\mathrm{d}\gamma^2;\,\,
\mathfrak{f}\in\mathcal{S}_{f,E}\right\} = \int_E
f\,\mathrm{d}\gamma^2,$$ which completes the proof. \qed \vskip
5pt

Similarly as in the Lebesgue measure and integral theory we
introduce sets which may be "neglected".

\begin{definition}\rm
Let $\gamma:\Sigma\to\Delta^+$ be a probabilistic-valued set
function. A set $E\in\Sigma$ is said to be \textit{$\gamma$-null},
if $\gamma_E=\varepsilon_0$. A set of all $\gamma$-null sets will
be denoted by $\mathcal{N}_\gamma$.
\end{definition}

Clearly, for each $\tau$-decomposable (sub)measure $\gamma$ w.r.t.
a triangle function $\tau$ the set $\mathcal{N}_\gamma$ is always
non-empty, because $\emptyset\in\mathcal{N}_\gamma$. Also, the
union of $\gamma$-null sets is a $\gamma$-null set. Indeed, if
$E_i$ are $\gamma$-null sets, then
$$\gamma_{\cup_{i=1}^{n}E_i}=\mathop{\bigoplus}\limits_{i=1}^{n}\gamma_{E_i}=
\mathop{\bigoplus}\limits_{i=1}^{n}\varepsilon_0=\varepsilon_0.$$
From antimonotonicity of a $\tau$-decomposable measure it follows
that each subset of $\gamma$-null set is again a $\gamma$-null
set.

\begin{example}\rm
Put $\gamma_E=\varepsilon_{\mu(E)}$. If $\Omega=\mathbb{N}$,
$\Sigma$ contains of all finite subsets of $\Omega$ and $\mu$ is a
counting measure on $\Sigma$, then the only $\gamma$-null set is
empty-set. On the other hand, for the Lebesgue measure $\mu$ on a
ring $\Sigma$ the set $\mathcal{N}_\gamma$ contains of all
Lebesgue $\mu$-measurable null sets from $\Sigma$.
\end{example}

\begin{example}\rm
Put $\gamma^\Phi_E=\mu(E)\odot \Phi$ for an arbitrary
$\Phi\in\Delta^+$. For $\Phi=\varepsilon_0$ all the sets from a
ring $\Sigma$ of subsets of $\Omega\neq\emptyset$ belong to
$\mathcal{N}_{\gamma^\Phi}$ for any finite measure $\mu$ on
$\Sigma$.
\end{example}

\begin{theorem}
Let $\tau\in\mathfrak{D}(\Delta^+)$ and $\gamma$ be a
$\tau$-decomposable measure on $\Sigma$. If
$E\in\mathcal{N}_\gamma$, then $\int_E f\,\mathrm{d}\gamma =
\varepsilon_0$ for any (simple) $\gamma$-integrable function
$f\in\Omega^{\overline{\mathbb{R}}_{+}}$.
\end{theorem}

\proof For a simple function of the form
$f=\sum\limits_{i=1}^{n}x_{i}\chi_{E_{i}}$ we have $$\int_E
f\,\mathrm{d}\gamma = \int_E
\left(\sum\limits_{i=1}^{n}x_{i}\chi_{E_{i}}\right)\,\mathrm{d}\gamma
= \mathop{\bigoplus}\limits_{i=1}^{n}x_i\odot\gamma_{E\cap
E_{i}}.$$ By antimonotonicity of $\gamma$ and $E\cap E_i\subseteq
E$ we get $\varepsilon_0\geq\gamma_{E\cap E_i}\geq
\gamma_E=\varepsilon_0$. Thus, $\gamma_{E\cap E_i}=\varepsilon_0$,
and then $$\int_E f\,\mathrm{d}\gamma =
\mathop{\bigoplus}\limits_{i=1}^{n}x_i\odot\gamma_{E\cap E_{i}} =
\mathop{\bigoplus}\limits_{i=1}^{n}x_i\odot\varepsilon_0 =
\varepsilon_0.$$

For a general non-negative $\gamma$-integrable function $f\in
\Omega^{\overline{\mathbb{R}}_{+}}$ there is a function
$H\in\Delta^+$ such that $\int_E
\mathfrak{f}\,\mathrm{d}\gamma\geq H$ for each
$\mathfrak{f}\in\mathcal{S}_{f,E}$ and $\int_E f\,\mathrm{d}\gamma
= \inf\{\int_E \mathfrak{f}\,\mathrm{d}\gamma;\,\,
\mathfrak{f}\in\mathcal{S}_{f,E}\}$. Since $\int_E
\mathfrak{f}\,\mathrm{d}\gamma = \varepsilon_0$ for each
$\mathfrak{f}\in\mathcal{S}_{f,E}$, then $\int_E
f\,\mathrm{d}\gamma = \varepsilon_0$. \qed \vskip 5 pt

\begin{corollary}\rm
Let $\tau\in\mathfrak{D}(\Delta^+)$ and $\gamma$ be a
$\tau$-decomposable measure on $\Sigma$. If
$\mathcal{N}_\gamma=\Sigma$, then $\int_E f\,\mathrm{d}\gamma =
\varepsilon_0$ for any (simple) $\gamma$-integrable function
$f\in\Omega^{\overline{\mathbb{R}}_{+}}$.
\end{corollary}


\section{Integral as a limit of integrals of simple functions}\label{sec_new}

In this section consider a $\sigma$-ring $\Sigma$ of subsets of
$\Omega\neq\emptyset$. Motivated by Halmos~\cite[\S 20, Theorem
B]{Halmos}, each non-negative (extended) real-valued function $f$
is the pointwise limit of a sequence $(f_n)_1^\infty$ of
non-decreasing simple functions, it is possible to define integral
via uniform approximation of a non-negative function by a sequence
of pointwise converging monotone sequence of simple functions.
However, we have to assume here that a triangle function $\tau$ is
continuous and a probabilistic-valued measure $\gamma$ is
continuous from below (or, it is sufficient to consider its
$\sigma$-additivity, see Lemma~\ref{lemma_continuity}). In what
follows $\mathfrak{CD}(\Delta^+)$ is the set of all continuous
distributive triangle functions on $\Delta^+$.

\begin{lemma}\label{lemma_simple}
Let $\tau\in \mathfrak{CD}(\Delta^+)$, $\gamma$ be a
$\tau$-decomposable measure on $\Sigma$ continuous from below and
$f\in\mathcal{S}$ . If $E_n\nearrow E$ , then
$$\lim\limits_{n\to\infty} \int_{E_n}f
\,\mathrm{d}\gamma=\int_{E} f\,\mathrm{d}\gamma.$$
\end{lemma}

\proof Let the simple function $f$ be of the form
$\sum\limits_{i=1}^m x_i\cdot\chi_{G_i}$, where
$x_i\in[0,+\infty)$, and $G_i\in\Sigma$ are pairwise disjoint
sets. From Lemma~\ref{lemma_continuity} and from continuity of
triangle function $\tau$ it follows
\begin{align*}
\lim\limits_{n\to\infty} \int_{E_n}f \,\mathrm{d}\gamma
&=\lim\limits_{n\to\infty} \int_{E_n}\left(\sum_{i=1}^m
x_i\cdot\chi_{G_i}\right)
\,\mathrm{d}\gamma=\lim\limits_{n\to\infty}\bigoplus\limits_{i=1}^m
x_i\odot\gamma_{E_n\cap G_i}=\bigoplus\limits_{i=1}^m
\lim\limits_{n\to\infty} x_i\odot\gamma_{E_n\cap G_i}\\& =
\bigoplus\limits_{i=1}^m \
x_i\odot\left(\lim\limits_{n\to\infty}\gamma_{E_n\cap G_i}\right)=
\bigoplus\limits_{i=1}^m \ x_i\odot\gamma_{E\cap
G_i}=\int_{E}f\,\mathrm{d}\gamma,
\end{align*}which completes the proof. \qed

Now we are able to prove that the $\gamma$-integral of a
measurable function $f$ may be written as limit of integrals of
simple functions $f_n$ pointwisely converging to $f$.

\begin{theorem}\label{thm_pointwise}
Let $\tau\in \mathfrak{CD}(\Delta^+)$, $\gamma$ be a
$\tau$-decomposable measure on $\Sigma$ continuous from below, and
$f\in\Omega^{\overline{\mathbb{R}}_+}$ be a $\gamma$-integrable
function  on a set $E\in\Sigma$. Then there exists a
non-decreasing sequence $(f_n)_1^\infty\in\mathcal{S}_{f,E}$
converging pointwisely to $f$ such that $$\int_E f
\mathrm{d}\gamma=\lim\limits_{n\to\infty}\int_E f_n
\mathrm{d}\gamma.$$
\end{theorem}

\proof By Halmos~\cite[\S 20, Theorem B]{Halmos} each non-negative
(extended) real-valued function $f$ is the pointwise limit of a
sequence $(f_n)_1^\infty$ of non-decreasing simple functions.
Since $\left(\int_E f_n\,\mathrm{d}\gamma\right)_1^\infty$ is a
non-increasing sequence of distance distribution functions, which
is bounded below by distance distribution function $\int_E f
\mathrm{d}\gamma$, existence of the limit follows. From
monotonicity of integral we immediately get the inequality
$\lim\limits_{n\to\infty}\int_E f_n \mathrm{d}\gamma\geq\int_E f
\mathrm{d}\gamma$. Thus, it suffices to show the reverse
inequality.

On the set $E\in\Sigma$ consider an arbitrary simple function
$\mathfrak{f}\in\mathcal{S}$ such that $\mathfrak{f}\leq f$. For a
fixed $t\in(0,1)$ put $E_n:=\{x\in E; f_n(x)\geq
t\cdot\mathfrak{f}(x)\}$. Clearly, $E_1\subseteq E_2\subseteq\dots
\subseteq E_n\subseteq\dots$ and
$\mathop{\cup}\limits_{n=1}^\infty E_n=E$. Thus,
$$\int_E f_n \mathrm{d}\gamma\leq\int_{E_n} f_n
\mathrm{d}\gamma\leq \int_{E_n} t\cdot\mathfrak{f}
\,\mathrm{d}\gamma=t\odot\int_{E_n}\mathfrak{f}
\,\mathrm{d}\gamma,$$ where the last equality is due to positive
homogeneity of the integral, see Theorem~\ref{nvlastnosti}. Then
Lemma~\ref{lemma_simple} yields $\lim\limits_{n\to\infty}
\int_{E_n}\mathfrak{f} \,\mathrm{d}\gamma =
\int_{E}\mathfrak{f}\,\mathrm{d}\gamma$. Therefore,
$$\lim\limits_{n\to\infty}\int_E f_n \,\mathrm{d}\gamma\leq
\lim\limits_{n\to\infty}t\odot\int_{E_n}\mathfrak{f}\,\mathrm{d}\gamma
=
t\odot\lim\limits_{n\to\infty}\int_{E_n}\mathfrak{f}\,\mathrm{d}\gamma=t\odot\int_{E}\mathfrak{f}\,\mathrm{d}\gamma.$$
Since $t\in(0,1)$, then $\lim\limits_{n\to\infty}\int_E
f_n\,\mathrm{d}\gamma\leq\int_{E}\mathfrak{f}\,\mathrm{d}\gamma$.
Taking the infimum of all such simple functions we get
$$\lim\limits_{n\to\infty}\int_E f_n\,\mathrm{d}\gamma\leq\inf\left\{\int_{E}\mathfrak{f}\,\mathrm{d}\gamma;\,\,\mathfrak{f}\in\mathcal{S}_{f,E}\right\}=\int_{E} f\,\mathrm{d}\gamma.$$\qed

The following theorem examines the linearity of the integral
\textbf{which depends on a choice of triangle function} $\tau$.

\begin{theorem}\label{vlastnosti}
Let $\tau\in \mathfrak{CD}(\Delta^+)$ and $\gamma$ be a
$\tau$-decomposable measure on $\Sigma$ continuous from below. If
$f,g\in \Omega^{\overline{\mathbb{R}}_{+}}$ are (simple)
$\gamma$-integrable functions on $E\in\Sigma$ and $\mathcal{I}:
f\mapsto\int_{\cdot}{f \,\mathrm{d}\gamma}$, then the
$\gamma$-integral is $\oplus_\tau$\,-linear, i.e.,
$\mathcal{I}(f+g)= \mathcal{I}(f)\oplus_\tau \mathcal{I}(g)$ if
and only if $\tau=\tau_{M}$. In general, $\mathcal{I}(f+g)\geq
\mathcal{I}(f)\oplus_\tau \mathcal{I}(g)$ if and only if
$\tau\geq\tau_{M}$.
\end{theorem}

\proof Let $f, g\in\mathcal{S}$. Then there exist pairwise
disjoint sets $B_1,\dots,B_p\in\Sigma$ such that both functions
may be expressed in the form $f=\sum\limits_{k=1}^p x_k
\chi_{B_k}$, $g=\sum\limits_{k=1}^p y_k \chi_{B_k}$. Thus, the
function $f+g$ is also a simple function of the form
$f+g=\sum\limits_{k=1}^p z_k\chi_{B_k}$ with $z_k=x_k+y_k$,
$k=1,2,\dots,p$.

(a) Let us consider $\tau=\tau_M$. Then
\begin{align*}
\int_E (f+g)\,\mathrm{d}\gamma & =  \bigoplus_{k=1}^p
(x_k+y_k)\odot\gamma_{E\cap B_k} = \left(\bigoplus_{k=1}^p
x_k\odot\gamma_{E\cap B_k}\right)\oplus_\tau
\left(\bigoplus_{k=1}^p y_k\odot\gamma_{E\cap B_k}\right) \\ & =
\int_E f\,\mathrm{d}\gamma \oplus_\tau \int_E g\,\mathrm{d}\gamma,
\end{align*}where the second equality holds if and only if $\tau=\tau_M$, see
Example~\ref{example MS}. \normalsize Now we show that $f+g$ is
$\gamma$-integrable on $E\in\Sigma$ provided that $f, g\in
\Omega^{\overline{\mathbb{R}}_{+}}$ are $\gamma$-integrable on
$E\in\Sigma$. By definition there exist $\varphi\in\Delta^+$ and
$\psi\in\Delta^+$ such that $\int_E \mathfrak{f}\,\mathrm{d}\gamma
\geq \varphi$ for each $\mathfrak{f}\in\mathcal{S}_{f,E}$ and
$\int_E \mathfrak{g}\,\mathrm{d}\gamma \geq \psi$ for each
$\mathfrak{g}\in\mathcal{S}_{g,E}$, respectively. Then
$\mathfrak{f}+\mathfrak{g}\in\mathcal{S}_{f+g,E}$ and
$$\int_E (\mathfrak{f}+\mathfrak{g})\,\mathrm{d}\gamma = \int_E
\mathfrak{f}\,\mathrm{d}\gamma \oplus_\tau \int_E
\mathfrak{g}\,\mathrm{d}\gamma \geq \varphi\oplus_\tau\psi$$ for
each $\mathfrak{f}+\mathfrak{g}\in\mathcal{S}_{f+g,E}$, i.e.,
$f+g$ is $\gamma$-integrable on $E\in\Sigma$. Since $\tau=\tau_M$
is continuous, because $M$ is a continuous t-norm, by
Theorem~\ref{thm_pointwise} there exist non-decreasing sequences
$(f_n)_1^\infty\in\mathcal{S}_{f,E}$ and
$(g_n)_1^\infty\in\mathcal{S}_{g,E}$ pointwisely converging to $f$
and $g$, respectively, such that
\begin{align*}
\int_E (f+g)\,\mathrm{d}\gamma & =  \lim_{n\to\infty} \int_E
(f_n+g_n)\,\mathrm{d}\gamma = \lim_{n\to\infty} \left(\int_E
f_n\,\mathrm{d}\gamma \bigoplus \int_E
g_n\,\mathrm{d}\gamma\right) \\ & = \lim_{n\to\infty} \int_E
f_n\,\mathrm{d}\gamma \bigoplus \lim_{n\to\infty}\int_E
g_n\,\mathrm{d}\gamma = \int_{E} f \,\mathrm{d}\gamma \bigoplus
\int_{E} g \,\mathrm{d}\gamma.
\end{align*}

(b) If $\tau\geq\tau_M$, then for $f,g\in\mathcal{S}$ we get $$
\int_E (f+g)\,\mathrm{d}\gamma \geq \left(\bigoplus_{k=1}^p
x_k\odot\gamma_{E\cap B_k}\right)\oplus_\tau
\left(\bigoplus_{k=1}^p y_k\odot\gamma_{E\cap B_k}\right) = \int_E
f\,\mathrm{d}\gamma \oplus_\tau \int_E g\,\mathrm{d}\gamma,
$$ where $\oplus_\tau$-suplinearity of $\gamma$-integral of simple
functions (the first inequality) is equivalent to
$\tau\geq\tau_M$, see Example~\ref{example MS}. The rest of the
proof is similar to the part (a). \qed \vskip 5pt

We illustrate the $\oplus_{\tau_{M}}$-linearity property of the
integral considering a class of probabilistic measures from
$\Gamma_A$.

\begin{figure}
\begin{center}
\includegraphics{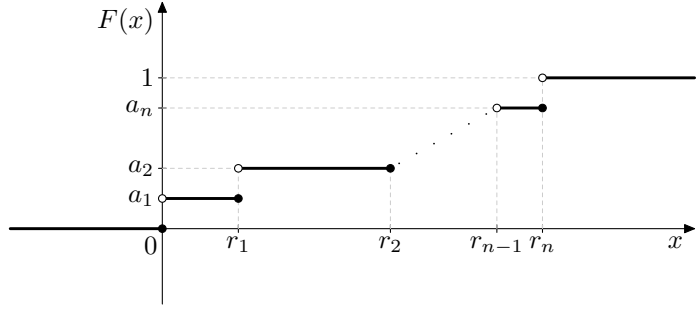}
\caption{Result of integration from Example~\ref{priklad2}}
\label{obr2}
\end{center}
\end{figure}

\begin{example}\rm\label{priklad2}
Let $\gamma^{a_i}\in\Gamma_A$ be probabilistic measures defined on
pairwise disjoint sets $E_i \in \Sigma$, where $i=1,2,\dots, n$
and $a_1\leq a_2\leq\dots\leq a_n$. Then the integral of function
$f= x_{0}\chi_{E}$ for $x_0\in \mathbb{R}_{+}$ (a constant
function on the set $E\in\Sigma$) computed on the sets $E_1,
E_2,\dots, E_n$ is a function given by
$$H=\mathop{\bigoplus}\limits_{i=1}^{n}\int_{E _i}f
\,\mathrm{d}\gamma^{a_i} = 
a_1\chi_{]0,r_1]}+a_2\chi_{]r_1,r_2]}+\dots
a_n\chi_{]r_{n-1},r_n]} + \chi_{]r_n,+\infty]},$$ where
$$r_i=s_1+\dots+s_i\quad \textrm{a} \quad s_i= \sum\limits_{k=1}^{i}x_{0}\mu(E\cap
E_{k})=\sum\limits_{k=1}^{i} \int_{E_k}f \mathrm{d}\mu.$$ Graph of
this function is illustrated in Fig.~\ref{obr2}.
\end{example}

In accordance with~\cite{Lpspaces} for a function $G\in\Delta^+$
and a constant $a\in[0,1]$ introduce the notation
$$G_{a}^{-}:=\inf\{x\in\mathbb{R};\, a\leq G(x)\},$$
$$G_{a}^{+}:=\sup\{x\in\mathbb{R};\, a\geq G(x)\}.$$ Then the distance
distribution function $H$ which is the result of integration in
Example~\ref{priklad2} may be written as
$$H=\mathop{\bigoplus}\limits_{i=1}^{n}\int_{E_i} f \,\mathrm{d}\gamma^{a_i}=
\sum\limits_{a\in A} a \chi_{]H_{a}^{-},H_{a}^{+}]}.$$ The
resulting form follows from observation that if $a=a_i$ for some
$i=1,2,\dots, n$, then $H_{a}^{-}<H_{a}^{+}$, otherwise
$H_{a}^{-}=H_{a}^{+}$.

\begin{theorem}
Let $\tau\in \mathfrak{CD}(\Delta^+)$ and $\gammaî$ be
$\tau$-decomposable measures on $\Sigma$ continuous from below
with $i\in\{1,2\}$. If $f\in \Omega^{\overline{\mathbb{R}}_{+}}$
is a (simple) $\gamma^i$-integrable function on $E\in\Sigma$ for
$i\in\{1,2\}$, then $f$ is
$\gamma^1\oplus_\tau\gamma^2$-integrable on $E\in\Sigma$, and it
holds
$$\int_E f\,\mathrm{d}(\gamma^1\oplus_\tau\gamma^2) = \int_E
f\,\mathrm{d}\gamma^1 \oplus_\tau \int_E f\,\mathrm{d}\gamma^2.$$
\end{theorem}

\proof First observe that if $\gamma^i$ are $\tau$-decomposable
measures on $\Sigma$ for $i\in\{1,2\}$, then by~\cite[Theorem
4.2]{HalHutMol} the set function
$\zeta:=\gamma^1\oplus_\tau\gamma^2$ is $\tau$-decomposable
measure on $\Sigma$ as well.

Let $f$ be a simple function of the form $f=\sum\limits_{i=1}^n
x_i\chi_{E_i}$. Then distributivity and associativity of $\tau$
yields
\begin{align*}
\int_E f\,\mathrm{d}\zeta & = \bigoplus_{i=1}^n
x_i\odot(\gamma^1_{E\cap E_i}\oplus_\tau\gamma^2_{E\cap E_i})
=\bigoplus_{i=1}^n (x_i\odot\gamma^1_{E\cap
E_i})\oplus_\tau(x_i\odot\gamma^2_{E\cap E_i}) \\ & =
\left(\bigoplus_{i=1}^n x_i\odot\gamma^1_{E\cap
E_i}\right)\oplus_\tau\left(\bigoplus_{i=1}^n
x_i\odot\gamma^2_{E\cap E_i}\right)= \int_E f\,\mathrm{d}\gamma^1
\oplus_\tau \int_E f\,\mathrm{d}\gamma^2.
\end{align*} Let us show that it holds for a non-negative integrable function
$f$. By Theorem~\ref{thm_pointwise} there is a non-decreasing
sequence $(f_n)_1^\infty\in\mathcal{S}_{f,E}$ pointwisely
converging to $f$ such that
\begin{align*}
\int_E f\,\mathrm{d}\zeta & =  \lim_{n\to\infty} \left(\int_E
f_n\,\mathrm{d}\gamma^1 \bigoplus \int_E
f_n\,\mathrm{d}\gamma^2\right) = \lim_{n\to\infty} \int_E
f_n\,\mathrm{d}\gamma^1 \bigoplus \lim_{n\to\infty} \int_E
f_n\,\mathrm{d}\gamma^2 \\ & =\int_E f\,\mathrm{d}\gamma^1
\bigoplus \int_E f\,\mathrm{d}\gamma^2,
\end{align*}where the continuity of $\oplus_\tau$ has been used.
\qed \vskip 5pt

As in the case of the Lebesgue integral, the $\gamma$-integral can
be seen as an extension of a $\tau$-decomposable measure.

\begin{theorem}\label{gammanadintegral}
Let $\tau\in\mathfrak{CD}(\Delta^+)$ and $\gamma$ be a
$\tau$-decomposable measure on $\Sigma$ continuous from below. For
each (simple) $\gamma$-integrable function
$f\in\Omega^{\overline{\mathbb{R}}_{+}}$ the set function $\nu^f:
\Sigma\to\Delta^+$ defined by $$\nu^f_{E} := \int_{E}{f
\,\mathrm{d}\gamma},\quad E\in\Sigma,$$ is a $\tau$-decomposable
measure.
\end{theorem}

\proof Clearly, $\nu^f_\emptyset=\varepsilon_0$. Let
$E,F\in\Sigma$ be disjoint sets. For a simple function
$f=\sum\limits_{i=1}^n x_i\chi_{E_i}$ we write
\begin{align*}
\nu^f_{E\cup F}&=\int_{E\cup F} f\,\mathrm{d}\gamma =
\bigoplus_{i=1}^n x_i\odot\gamma_{(E\cup F)\cap E_i} =
\bigoplus_{i=1}^n x_i\odot\gamma_{(E\cap E_i)\cup(F\cap E_i)} \\
&= \bigoplus_{i=1}^n x_i\odot(\gamma_{E\cap
E_i}\oplus_\tau\gamma_{F\cap E_i}) = \left(\bigoplus_{i=1}^n
x_i\odot\gamma_{E\cap
E_i}\right)\oplus_\tau\left(\bigoplus_{i=1}^n
x_i\odot\gamma_{F\cap E_i}\right) \\ & = \nu^f_{E}\oplus_\tau
\nu^f_{F}.
\end{align*}For a general $\gamma$-integrable function $f$ by Theorem~\ref{thm_pointwise} there exists a non-decreasing sequence $(f_n)_1^\infty\in\mathcal{S}_{f,E\cup F}$ such that
\begin{align*}
\nu^{f}_{E\cup F} & = \lim_{n\to\infty}\nu^{f_n}_{E\cup F} =
\lim_{n\to\infty} \left(\nu^{f_n}_E \oplus_\tau \nu^{f_n}_F\right)
= \left(\lim_{n\to\infty}\nu^{f_n}_E\right) \oplus_\tau
\left(\lim_{n\to\infty}\nu^{f_n}_F\right) = \nu^f_{E}\oplus_\tau
\nu^f_{F},
\end{align*}which completes the proof. \qed \vskip 5pt

\begin{remark}\rm
According to Proposition~\ref{proposition1}(ii) the
$\gamma$-integral is antimonotone, i.e., the inequality
$\nu^f_E\geq\nu^f_F$ holds for each (simple) $\gamma$-integrable
function $f$ and each $E,F\in\Sigma$ such that $E\subseteq F$.
Moreover, the set function $\nu^f$ and the operator $\mathcal{I}$
are related by the formula $\nu^f_E = \mathcal{I}(f\cdot\chi_E)$.
Theorem~\ref{gammanadintegral} shows that the $\gamma$-integral is
indeed a proper extension of the underlying measure $\gamma$
because of $\nu^{\chi_E} = \gamma_E$ for each $E\in\Sigma$.
\end{remark}

\begin{problem}
Under which conditions on a probabilistic-valued set function
$\gamma$, triangle function $\tau$ and functions
$f,g\in\Omega^{\overline{\mathbb{R}}_+}$ the following equalities
hold
$$\int_E f\,\mathrm{d}\nu^g = \int_E g\,\mathrm{d}\nu^f = \int_E f\cdot g\,\mathrm{d}\gamma, \quad E\in\Sigma?$$
\end{problem}

\begin{definition}\rm
Functions $f,g\in\Omega^{\overline{\mathbb{R}}_{+}}$ are said to
be \textit{equal a.e.} on a set $S$, we write $f=g$ a.e., if
$f(x)=g(x)$ for all $x\in S\setminus E$, where
$E\in\mathcal{N}_\gamma$.
\end{definition}

\begin{theorem}
Let $\tau\in\mathfrak{CD}(\Delta^+)$, $\gamma$ be a
$\tau$-decomposable measure on $\Sigma$ continuous from below, and
$f,g\in\Omega^{\overline{\mathbb{R}}_{+}}$ such that $f=g$ a.e. on
a set $E\in\Sigma$. Then $\int_E
f\,\mathrm{d}\gamma=\int_{E}g\,\mathrm{d}\gamma$.
\end{theorem}

\proof Let $F$ be a $\gamma$-null set and $f(x)=g(x)$ for all
$x\in E\setminus F$. Since by Theorem~\ref{gammanadintegral} the
$\gamma$-integral is a  $\tau$-decomposable measure, then
\begin{align*}\int_E f\,\mathrm{d}\gamma & = \int_{(E\setminus
F)\cup F}f\,\mathrm{d}\gamma = \int_{(E\setminus
F)}f\,\mathrm{d}\gamma\oplus_\tau\int_{F}f\,\mathrm{d}\gamma =
\int_{E\setminus F}f\,\mathrm{d}\gamma=\int_{E\setminus
F}g\,\mathrm{d}\gamma = \int_{E\setminus
F}g\,\mathrm{d}\gamma\oplus\varepsilon_0 \\ & = \int_{(E\setminus
F)}g\,\mathrm{d}\gamma\oplus\int_{F}g\,\mathrm{d}\gamma=\int_{E}g\,\mathrm{d}\gamma,\end{align*}
which completes the proof. \qed \vskip 5 pt

\section*{Concluding remarks}

Origin of probabilistic-valued set functions comes from the fact
that they work in such situations in which we have only a
\textit{probabilistic information} about measure of a set (recall
a similar situation in the framework of information measures as
discussed in~\cite{KdF74}). For example, in Moore's interval
analysis, if rounding of reals is considered, then the uniform
distributions over intervals describe our information about the
measure of a set. Another such probabilistic information occurs in
biometric decision-making where the information is often obtained
from biometric sensors as well as from the analysis of historical
data. In this paper we have discussed an approach to the
investigation of a Lebesgue-type integral w.r.t.
probabilistic-valued set functions. Although the introduced
probabilistic integral is based on a $\tau$-decomposable measure,
all the results may easily be rewritten for a $\tau$-decomposable
submeasure leading to a certain non-additive integral. Some
further research of the integral will be concentrated on the
problems of convergences for sequences of measurable functions,
convergences of the corresponding probabilistic-valued integrals
as well as their possible applications in different mathematical
theories.

The applicability of our integral in approximate reasoning will be
illustrated on a simple example. Consider a universe $\Omega =
\{1,\dots,n\}$ of criteria (and then $\Sigma = 2^\Omega$ is the
power set of $\Omega$). Consider a fixed distributive triangle
function $\tau$. To define properly a probabilistic-valued
$\oplus$-decomposable measure (where $\oplus$ is defined by
$\tau$), it is enough to have the knowledge about measures of
singletons of $\Omega$. For each such singleton $E_i=\{i\}$, let
$\gamma_{E_i}$ be a distance distribution function of a random
variable uniformly distributed over non-negative interval
$[a_i,b_i]$, i.e., $\gamma_{E_i} = F_{[a_i,b_i]}$. For any
alternative $\textbf{a}$, its score vector $\textbf{x}$ in
$[0,1]^n$ is assigned. Now, we can apply our integral to determine
a probabilistic-valued utility function,
$$U(\textbf{a}) = \int \textbf{x}\,\mathrm{d}\gamma.$$
For the sake of simplicity, we illustrate only one distinguished
case: let $T = M$ be the strongest minimum t-norm, and let $\tau =
\tau_{M}$. Then $\tau$ is a continuous distributive distance
distribution function, and $$U(\textbf{a}) = \int
\textbf{x}\,\mathrm{d}\gamma = F_{[\alpha, \beta]},$$ where
$\alpha = \sum x_i\cdot a_i$ and $\beta = \sum x_i\cdot b_i$. As
we can see, in this case we have recovered standard Moore's
weighted arithmetic mean with interval-valued means.

It is worth to be mentioned that a different approach to
probabilistic-integral may be provided, e.g. via the Choquet one.
For a partition $D=\{a_0, a_1, \dots, a_n\}$ with
$0=a_0<a_1<\dots<a_n<+\infty$ we define a corresponding simple
function $f_D: \Omega\to[0,+\infty[$ as follows $$f_D=\sum_{i=1}^n
(a_i-a_{i-1})\chi_{E_i},$$ where $E_i=\{x\in\Omega; f_D(x)\geq
a_i\}$. Then a Choquet-like integral of a simple function $f_D$ on
a set $E\in\Sigma$ w.r.t. a $\tau$-decomposable measure is defined
by the equality
$$(C)\int_E f_D\,\mathrm{d}\gamma := \bigoplus_{i=1}^n (a_i-a_{i-1})\odot\gamma_{E\cap E_i}.$$
The integral of a non-negative measurable function $f$ may be
defined as $$(C)\int_E f\,\mathrm{d}\gamma =
\lim_{D\in\mathcal{D}} (C)\int_E f_D\,\mathrm{d}\gamma$$ provided
that the limit (in the sense of Moore-Smith) exists. Here,
$\mathcal{D}$ is the set of all partitions of $\mathbb{R}_+$ with
the order induced by inclusion. Note that this approach may be
understood as a modification of Aumann integral, or, of
Choquet-type integral based on interval-valued measures as
discussed e.g. in~\cite{Jang}.

\section*{Acknowledgement}

We thank Radko Mesiar for his insightful comments and remarks on
Moore's interval analysis. We acknowledge a partial support of
grants VEGA 1/0171/12, VVGS-PF-2013-115 and VVGS-2013-121.


\vspace{5mm}

\noindent \small{\textsc{Lenka Hal\v cinov\'a, Ondrej Hutn\'ik
} \newline Institute of Mathematics, Faculty of Science, Pavol
Jozef \v Saf\'arik University in Ko\v sice,
\newline {\it Current address:} Jesenn\'a 5, SK 040~01 Ko\v sice,
Slovakia,
\newline {\it E-mail addresses:} lenka.halcinova@student.upjs.sk
\newline \phantom{{\it E-mail addresses:}} ondrej.hutnik@upjs.sk

\end{document}